\documentclass[11pt]{article}
\usepackage[final]{epsfig}
\usepackage{graphics}
\usepackage{amsmath}
\usepackage{amsfonts}
\usepackage{latexsym}
\usepackage{amssymb}
\usepackage{graphicx}
\usepackage{url}
\usepackage{epstopdf}

\begin{document}
\newcommand{\eps}{{\varepsilon}}
\newcommand{\proofend}{$\Box$\bigskip}
\newcommand{\C}{{\mathbb C}}
\newcommand{\Q}{{\mathbb Q}}
\newcommand{\R}{{\mathbb R}}
\newcommand{\Z}{{\mathbb Z}}
\newcommand{\RP}{{\mathbb {RP}}}
\newcommand{\CP}{{\mathbb {CP}}}
\newcommand{\Tr}{\rm Tr}
\def\proof{\paragraph{Proof.}}

\newcommand{\SL}{\operatorname{SL}}
\newcommand{\PGL}{\operatorname{PSL}}

\title{A four vertex theorem for frieze patterns?}

\author{Serge Tabachnikov\footnote{
Department of Mathematics,
Penn State University,
University Park, PA 16802;
tabachni@math.psu.edu}
}

\date{}
\maketitle

\section{The four vertex theorem}

\begin{figure}[hbtp]
\centering
\includegraphics[height=2in]{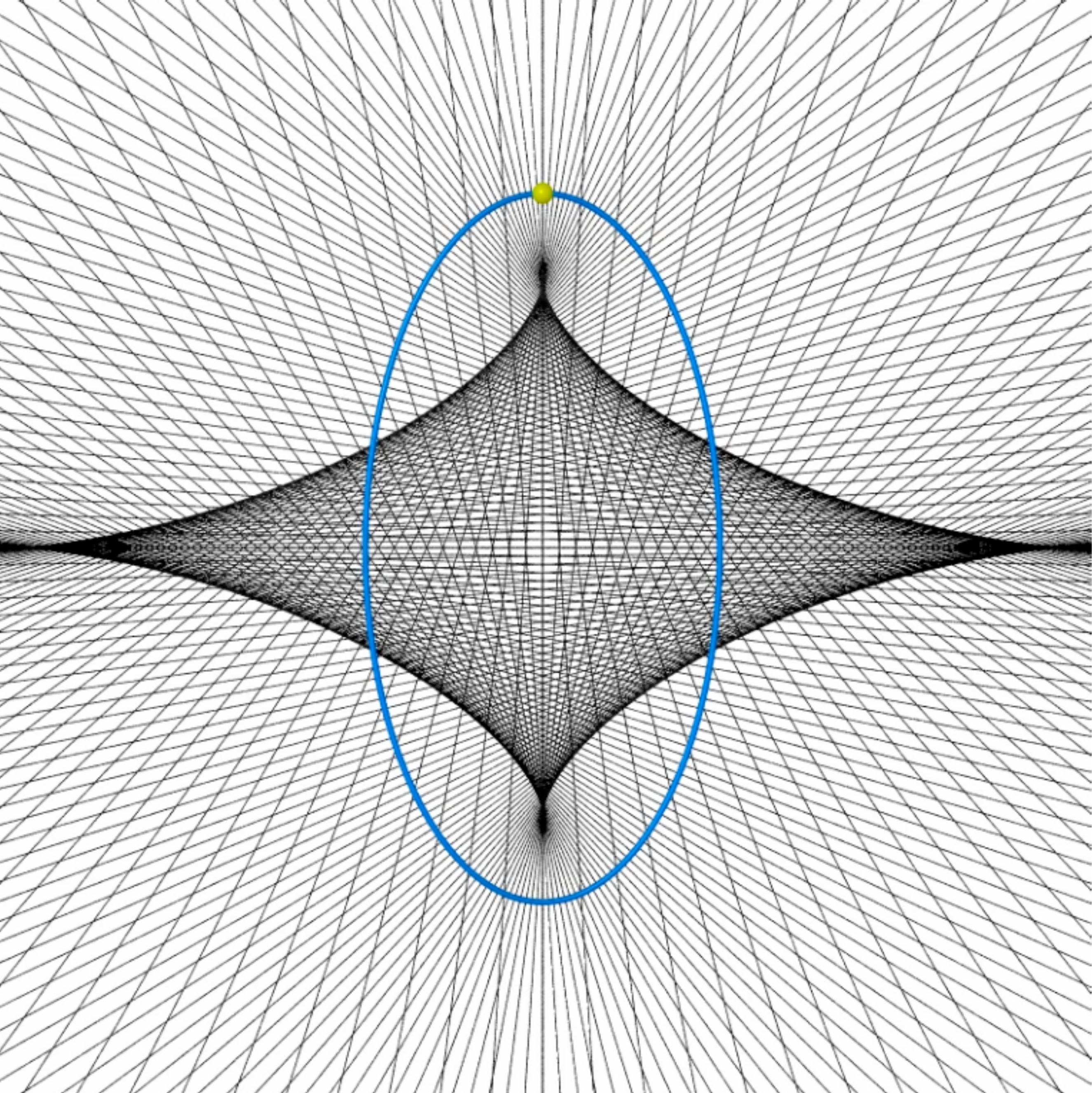}
\caption{An ellipse and its evolute, the envelope of its normals or, equivalently, the locus of the centers of curvature. The cusps of the evolute correspond to the vertices of the curve.}
\label{ellipse}
\end{figure}

The classic 4-vertex theorem states that {\it the curvature of a smooth closed convex planar curve has at least four critical points}, see Figure \ref{ellipse} for an illustration.

Since its discovery by  S. Mukhopadhyaya in 1909, this theorem has generated a large literature, comprising various generalizations and variants of this result; see  \cite{DGPV,GTT,OT3,Pak} for a sampler. 

One approach to the proof of the 4-vertex theorem is based on the following observation: if a $2\pi$-periodic function $f(x)$ is $L^2$-orthogonal to the first harmonics, that is, to the functions $1, \sin x, \cos x$, then $f(x)$ must have at least four sign changes over the period. 

The proof is simple: since $\int_0^{2\pi} f(x) dx =0$, the function $f(x)$ must change sign. If there are only two sign changes, one can find a linear combination $g(x)=c+a\cos x + b \sin x$ that changes sign at the same points as $f(x)$. Since the first harmonic $g(x)$ cannot have more than two sign changes, $f(x) g(x)$ has a constant sign, and $\int_0^{2\pi} f(x) g(x) dx \neq 0$, a contradiction. Discrete versions of this argument are in the hearts of our proofs presented below.

(In the 4-vertex theorem, one takes $f(x)=p'(x)+p'''(x)$, where $p(x)$ is the support function of the curve; then $p(x)+p''(x)$ is the curvature radius.) 

The above observation is a particular case of the Sturm-Hurwitz theorem: {\it the number of zeros of a  periodic function is not less than the number of zeros of its first non-trivial harmonic}, see \cite{OT3} for five proofs and applications of this remarkable result.

\begin{figure}[hbtp]
\centering
\includegraphics[height=2in]{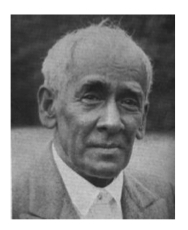}\quad
\includegraphics[height=2in]{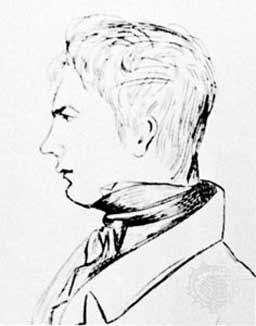}\quad
\includegraphics[height=2in]{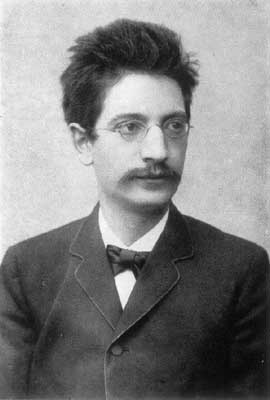}
\caption{Syamadas Mukhopadhyaya, Jacques Charles Fran\c{c}ois Sturm, and Adolf Hurwitz.}
\label{MSH}
\end{figure}

\section{Frieze patterns}

A frieze pattern is an array of numbers consisting  of finitely many bi-infinite rows; each next row is offset half-step from the previous one. The top two rows consist of 0s and of 1s, respectively, the bottom two rows are the row of 1s and 0s as well, and every elementary diamond
$$
 \begin{matrix}
 &N&
 \\
W&&E
 \\
&S&
\end{matrix}
$$
satisfies the unimodular relation $EW-NS=1$. The number of non-trivial rows is called the width of a frieze pattern. Denote the width by $w$ and set $n=w+3$.

\begin{figure}[hbtp]
\centering
\includegraphics[height=2.5in]{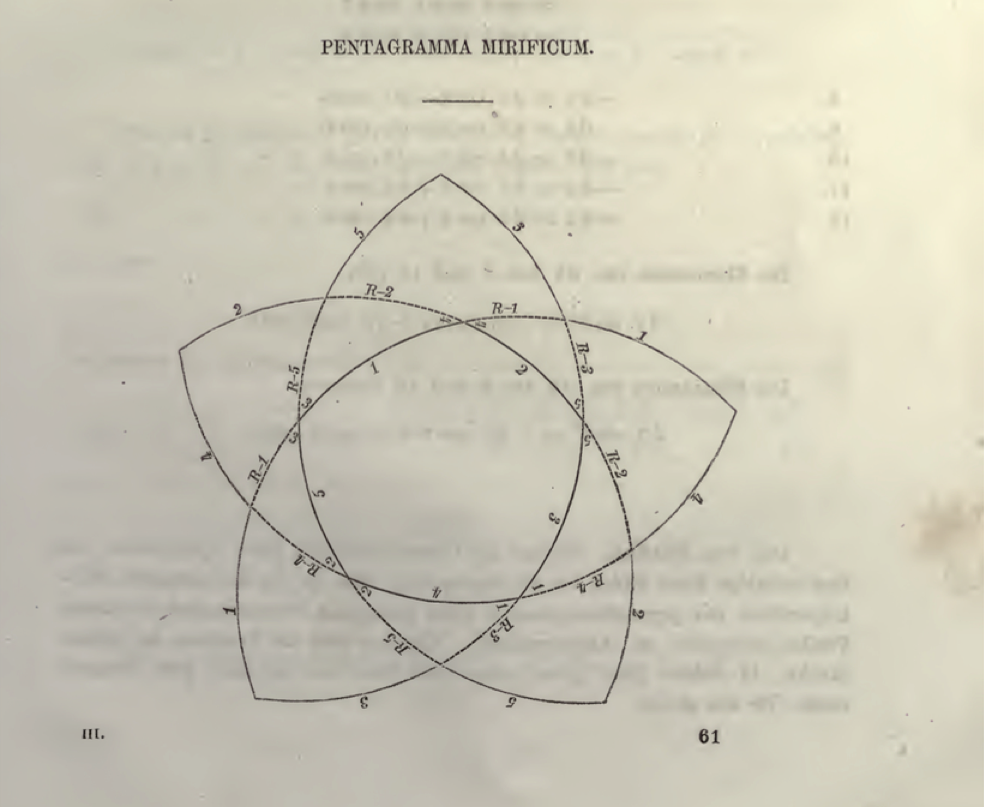}
\caption{Pentagramma mirificum of Carl Friedrich Gauss.}
\label{miri}
\end{figure}

For example,  a general frieze pattern with $w=2, n=5$ looks like this:
$$
 \begin{array}{ccccccccccc}
\cdots&&1&& 1&&1&&\cdots
 \\[4pt]
&x_1&&\frac{x_2+1}{x_1}&&\frac{x_1+1}{x_2}&&x_2&&
 \\[4pt]
\cdots&&x_2&&\frac{x_1+x_2+1}{x_1x_2}&&x_1&&\cdots
 \\[4pt]
&1&&1&&1&&1&&
\end{array}
$$
where the rows of 0s are omitted. These formulas appeared in the paper by Gauss ``Pentagramma Mirificum", published posthumously; Gauss calculated  geometric quantities characterizing spherical self-polar pentagons, see Figure \ref{miri}. See also A. Cayley's paper  \cite{Cay}. (According to Coxeter \cite{Cox} -- the very paper where frieze patterns were introduced -- the story goes further back, to N. Torporley, who in 1602 investigated the five ``parts" of a right-angled spherical triangle, anticipating by a dozen years the rule of J. Napier in spherical trigonometry.) 

And here is a frieze pattern of width four whose entries are natural numbers:
$$
 \begin{array}{ccccccccccccccccccccccc}
&&1&&1&&1&&1&&1&&1&&1&&1
 \\[4pt]
&1&&3&&2&&2&&1&&4&&2&&1&
 \\[4pt]
&&2&&5&&3&&1&&3&&7&&1&&2&
 \\[4pt]
&1&&3&&7&&1&&2&&5&&3&&1&
\\[4pt]
&&1&&4&&2&&1&&3&&2&&2&&1&
\\[4pt]
&1&&1&&1&&1&&1&&1&&1&&1
\end{array}
$$
The very existence of such frieze patterns is surprising: the unimodular rule $EW-NS=1$ does not agree easily with the property of being a positive integer! 

\begin{figure}[hbtp]
\centering
\includegraphics[height=2in]{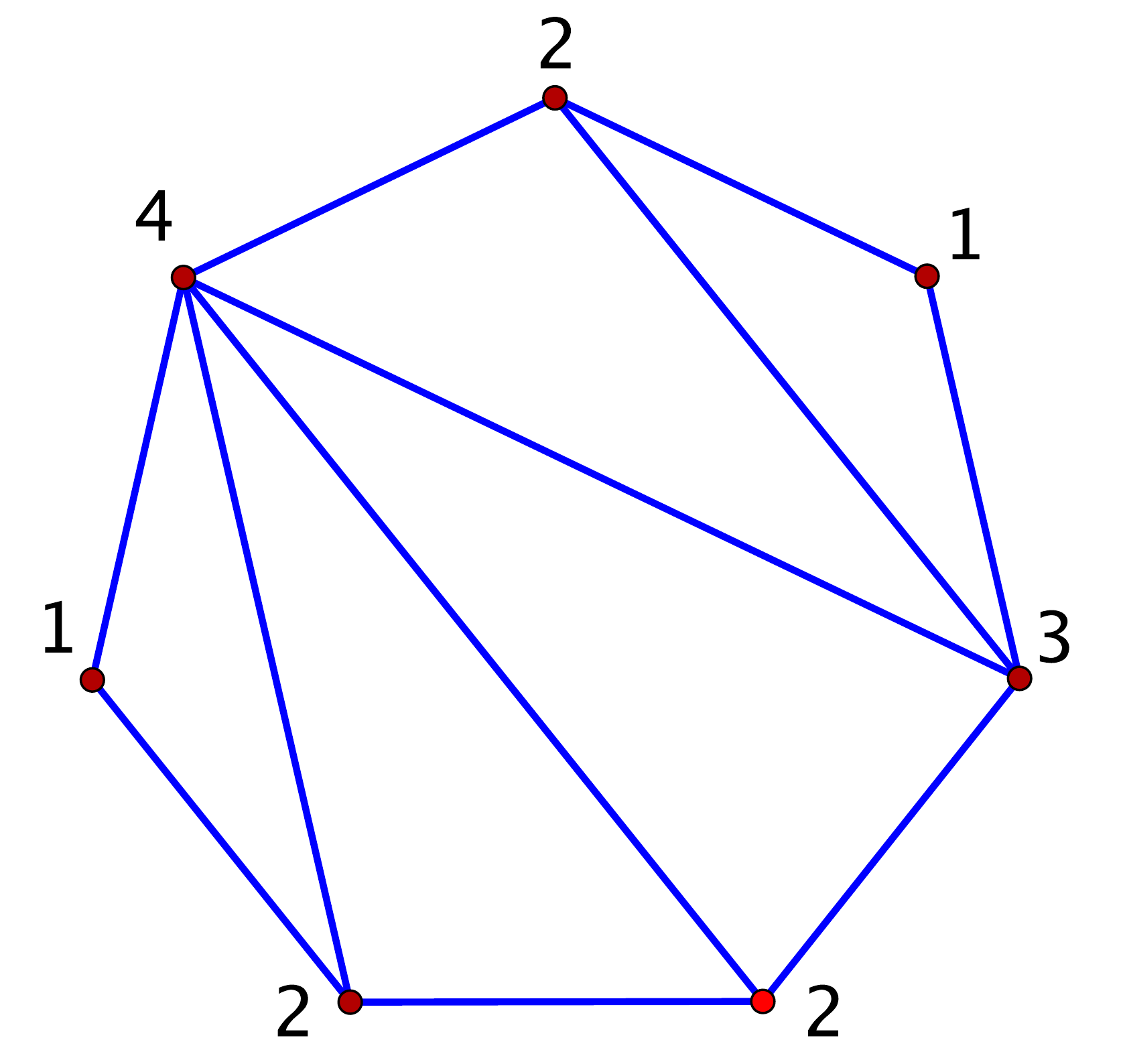}
\caption{A triangulation of a heptagon: the labels are the number of the triangles adjacent to each vertex. These numbers comprise the first row of the frieze pattern.}
\label{heptagon}
\end{figure}

The  frieze patters consisting of positive integers were classified by Conway and Coxeter \cite{CoCo}: they are in 1-1 correspondence with the triangulations of a convex $n$-gons by diagonals, and there are
$\frac{(2(w+1))!}{(w+1)!(w+2)!}$ (Catalan number) of them;  see \cite{Hen} for an exposition of this beautiful theorem. For example, the above frieze pattern corresponds to the triangulation in Figure \ref{heptagon}.

For a while, frieze patterns remained a relatively esoteric subject, but recently they have attracted much attention due of their significance in algebraic combinatorics and the theory of cluster algebras. I recommend a  comprehensive contemporary survey of this subject \cite{Mor}. 

\begin{figure}[hbtp]
\centering
\includegraphics[height=1.5in]{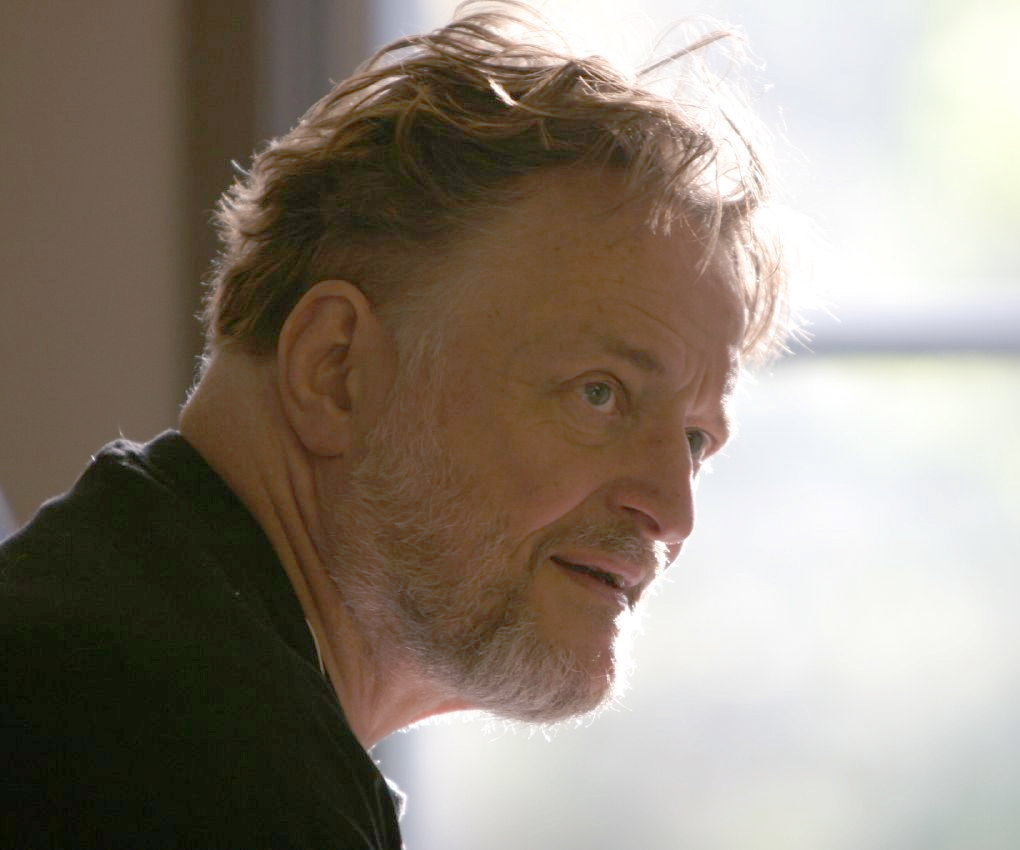}\quad
\includegraphics[height=1.5in]{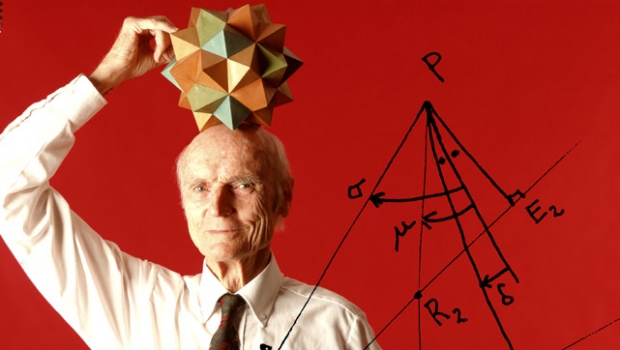}
\caption{John Horton Conway and Harold Scott MacDonald Coxeter.}
\label{ConCox}
\end{figure}

Let us summarize the basic properties of frieze patterns relevant to this article. Denote by $a_i$  the entries of the first non-trivial row.
\begin{enumerate}
\item The NE diagonals of a frieze pattern satisfy the 2nd order linear recurrence (discrete Hill's equation)
$$V_{i+1} = a_i V_i - V_{i-1}$$
 with $n$-periodic coefficients whose all solutions are antiperiodic, i.e., $V_{i+n}=-V_i$ for all $i$:
 $$
\begin{array}{ccccccccccc}
0\ \ &&0&&0&&0\\
&1&&1&&1\\
&&a_1&&a_{2}&&a_{3}\\
&&&a_1a_2-1&&a_2a_3-1\\
&&\cdots&&a_1a_2a_3-a_1-a_3&&\cdots\\
\end{array}
$$
\item The solutions of the discrete Hill's equation can be thought of as  polygonal lines $\ldots,V_1, V_2, \ldots$ in $\R^2$, with $\det (V_i,V_{i+1})=1$ and  $V_{i+n}=-V_i$. Such polygonal line is well defined up to $\SL(2,\R)$-action. The projections of the vectors $V_i$ to $\RP^1$ form an $n$-gon therein, well-defined up to a M\"obius transformation. For odd $n$, this correspondence between frieze patterns of width $n-3$ and projective equivalence classes of $n$-gons in the projective line is 1-1. 
\item Label the entries as follows:
$$
\begin{array}{ccc}
&v_{i,j}&\\
v_{i,j-1}&&v_{i+1,j}\\
&v_{i+1,j-1}&
\end{array}
$$
with $a_i=v_{i-1,i+1}$. 
Then one has $v_{i,j}=\det(V_i,V_j)$, explaining the glide reflection symmetry of the entries: $v_{i,j} = v_{j,i+n}$ 

The Conway-Coxeter article \cite{CoCo} starts with a description of the seven ornamental frieze patterns  where the glide reflection symmetry is represented by 
$
{\bf \large
\ldots b\quad p\quad b\quad p\quad b\quad p\ldots
}
$
and described as ``the relation between successive footprints when one walks along a straight path covered with snow". In Conway's nomenclature, this ornamental frieze pattern is called ``step", see Figure \ref{step}.

\begin{figure}[hbtp]
\centering
\includegraphics[height=.7in]{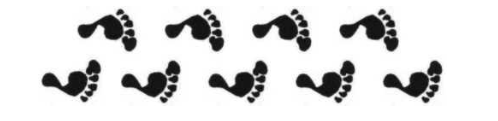}
\caption{An ornamental frieze pattern with the glide reflection symmetry.}
\label{step}
\end{figure}

\item The entries of a frieze pattern are given by the 3-diagonal determinants 
$$
\det\left|\begin{array}{cccccc}
a_{j}&1&&&\\\
1&a_{j+1}&1&&\\
&\ddots&\ddots&\ddots&\\
&&1&a_{i-1}&1\\
&&&1&a_{i}
\end{array}\right|,
$$
the continuants (called so because of their relation with continued fractions; see \cite{CO} for an intriguing history of this name).
\end{enumerate}

\section{A problem, a theorem, and a counter-example}

I shall be concerned with frieze patterns whose entries are positive real numbers. Given two such frieze patterns of the same width $w$,  choose a row and consider the $n$-periodic sequence of the differences of  the respective entries of the two friezes. I am interested in the number of sign changes in this sequence over the period. 

More precisely, let $1\le k \le [w/2]$ be the number of a row (we don't need to go beyond $[w/2]$ due to the glide symmetry), and let 
$v_{i,i+k+1}$ and $u_{i,i+k+1}$ be the entries of $k$th rows of the two frieze patterns. I am interested in the sign changes of $v_{i,i+k+1}-u_{i,i+k+1}$ as $i$ increases by 1 (not excluding the case when either of these differences vanishes).
\medskip

\noindent {\bf Problem 1.}
{\it 
For which $k$ must the cyclic sequence $v_{i,i+k+1}-u_{i,i+k+1}$ have at least four sign changes?
}
\medskip

As a partial answer, one has
\medskip

\noindent {\bf Theorem.}
{\it 
Four sign changes must occur for $k=1$ and for $k=2$.  In addition, for every $k$, four sign changes must occur in the 
infinitesimal version of the problem.
}
\medskip

Let me explain the last statement. 

Consider a frieze pattern whose first row is constant: $a_i=2x$ for all $i$. Then each next row is constant as well, and their entries, denoted by $U_k(x)$, satisfy $U_{k+1} = 2x U_k (x)- U_{k-1}(x)$ with $U_0(x)=1, U_1(x)=2x$. That is, $U_k(x)$ is the Chebyshev polynomial of the second kind:
$$
U_k(\cos \alpha) = \frac{\sin(k+1)\alpha}{\sin \alpha}.
$$
For this constant frieze pattern to have width $n-3$, set $\alpha=\pi/n$.

For the infinitesimal version of Problem 1,  take this constant frieze pattern and its infinitesimal deformation in the class of frieze patterns. 

Originally, I hoped that Problem 1 had an affirmative answer for all values of $k$. However, this conjecture was over-optimistic. The following counter-example is provided by Michael Cuntz; in this example, $w=5$ (the smallest possible not to contradict Theorem), all entries are positive rational numbers, and the differences of the entries of  the third row are all positive (this row is 4-periodic due to the glide symmetry). I present only the first lines of the two frieze patterns; these are 8-periodic sequences:
$$
\left(2,\ 2,\ 4,\ 2,\ 3,\ \frac{18}{41},\ 41, \frac{30}{41}\right) \ \ {\rm and}\ \ 
\left( 5,\ \frac{21}{97},\ 194,\ \frac{36}{97},\ 3,\ 5,\ 1,\ 5 \right).
$$
It still may be  possible that the bold conjecture holds for Conway-Coxeter frieze patterns that consist of positive integers.

\section{Proofs}

\paragraph{Case $k=1$.} 
 Let $a_i$ and $b_i$ be the entries of the first rows of the two frieze patterns. Consider the respective discrete Hill's equations
$$
V_{i+1} = a_i V_i - V_{i-1},\  U_{i+1} = b_i U_i - U_{i-1}.
$$
Let $U_i$ and $V_i$ be some solutions.I claim that the sequence $a_i-b_i$ is $\ell_2$-orthogonal to $U_i V_i$:
\begin{equation} \label{orth}
\sum _1^n (a_i - b_i) U_i V_i = 0. 
\end{equation}
Indeed, 
$$
\sum _1^n (a_i - b_i) U_i V_i = \sum _1^n [U_{i+1} + U_{i-1}] V_i - U_i [V_{i+1} + V_{i-1}] = 0,
$$
due to antiperiodicity.

Note that the space of solutions of a discrete Hill equation is 2-dimensional, and that its solutions are non-oscillating  in the sense that they change sign only once over the period (since the entries of the frieze pattern are positive). 

Assume that $a_i - b_i$ does not change sign at all. Choose the initial conditions for solutions $U_{i}$ and $V_i$ as follows:
$$
U_1=-1,U_2=1,V_1=-1,V_2=1.
$$
That is, both solution change sign from $i=1$ to $i=2$, and then, due to the non-oscillating property, there are no other sign changes. Hence $U_i V_i >0$ for all $i$, contradicting (\ref{orth}).

Likewise, if $a_i - b_i$ changes sign only twice, from $i_1$ to $i_1+1$, and from $i_2$ to $i_2+1$, choose the initial conditions for solutions $U_i$ and $V_i$ as follows:
$$
U_{i_1}=-1,U_{i_1+1}=1,V_{i_2}=-1,V_{i_2+1}=1.
$$
Then $(a_i - b_i) U_i V_i$ has a constant sign, again contradicting  (\ref{orth}).
\proofend

This result, along with its proof, is a discrete version of the following theorem from \cite{OT1} concerning Hill's equations $\varphi''(x) = k(x) \varphi(x)$ whose solutions are $\pi$-antiperiodic (and hence the potential $k(x)$ is $\pi$-periodic) and disconjugate, meaning that every solution changes sign only once on the period $[0,\pi)$. The claim is that, {\it given two such equations, the function $k_1(x) - k_2(x)$ has at least four zeroes on $[0,\pi)$}.

This theorem is equivalent to the beautiful theorem of E. Ghys: {\it the Schwartzian derivative of a diffeomorphism of $\RP^1$ has at least four distinct zeroes}, see \cite{OT3} for the relation of the Schwartzian derivative with the Hill equation, and an explanation why zeroes of the Schwartzian derivative are the vertices of a curve in  Lorentzian geometry. 

\paragraph{Case $k=2$.} 
As I mentioned, to a frieze pattern there corresponds an $n$-gon in $\RP^1$. The entries of the second row of the frieze pattern are the cross-ratios of the consecutive quadruples of the vertices of this $n$-gon, where cross-ratio is defined as
$$
[a,b,c,d]_1 = \frac{(d-a)(c-b)}{(d-c)(b-a)},
$$
see \cite{Mor}.

On the other hand, one of the results in \cite{OT2}, another discretization of Ghys's theorem, states that, given two cyclically ordered  $n$-tuples of points $x_i$ and $y_i$ in $\RP^1$, the difference of the cross-ratios $[x_i,x_{i+1},x_{i+2},x_{i+3}]_2  - [y_i,y_{i+1},y_{i+2},y_{i+3}]_2$ changes sign at least four times; here the cross-ratio is defined by
$$
[a,b,c,d]_2 = \frac{(d-b)(c-a)}{(d-c)(b-a)}.
$$
To complete the proof, observe that $[a,b,c,d]_2 - [a,b,c,d]_1 =1$.
\proofend

\paragraph{Infinitesimal version, $k$ arbitrary.}
Consider the polygonal line 
$$
V_i = \frac{1}{\sqrt{\sin \frac{\pi}{n}}} \left(\cos \frac{i\pi}{n}, \sin \frac{i\pi}{n}  \right),
$$
so that $[V_i,V_{i+1}]=1$ and $V_{i+n}=-V_i$ hold.

Let 
$$
W_i = V_i + \eps E_i,\ [W_i,W_{i+1}]=1
$$
be an infinitesimal deformation of this polygon $V_i$. I assume in our calculations that $\eps^2=0$.

Let
\begin{equation} \label{EtoV}
E_i=p_iV_i+\bar p_i V_{i+1} = q_i V_i + \bar q_i V_{i-1}.
\end{equation}
We shall express the $n$-periodic coefficients $p_i,\bar p_i, \bar q_i$ via the coefficients $q_i$, that solely determine the deformation.

To do so, use the fact that $V_{i+1}=cV_i-V_{i-1}$ with $c= 2 \cos (\pi/{n})$.
This linear relation must be equivalent to the second equality in (\ref{EtoV}), hence
$$
q_i-p_i=c \bar p_i,\ \bar p_i = - \bar q_i.
$$
We also have $[W_i,W_{i+1}]=1$,  implying $[V_i,E_{i+1}]+[E_i,V_{i+1}]=0$ and, using (\ref{EtoV}),
$p_i=-q_{i+1}$. Thus
\begin{equation} \label{viaq}
p_i = -q_{i+1},\ \bar p_i = \frac{1}{c} (q_i+q_{i+1}),\ \bar q_i = -\frac{1}{c} (q_i+q_{i+1}).
\end{equation}

Now fix $k$ and consider the deformation of  $(k-1)$st row of the frieze pattern:
$$
[W_i,W_{i+k}]=[V_i,V_{i+k}] + \eps ([V_i,E_{i+k}]+[E_i,V_{i+k}]).
$$ 
Using (\ref{EtoV}) and (\ref{viaq}), one finds
\begin{equation*}
\begin{split}
[V_i,E_{i+k}]+[E_i,V_{i+k}] = (q_{i+k}-q_{i+1}) \left([V_i,V_{i+k}] - \frac{1}{c} [V_i,V_{i+k-1}]\right)&\\
 - \frac{1}{c} (q_{i+k+1}-q_i) [V_i,V_{i+k-1}]&\\
= \frac{1}{\sin \frac{2\pi}{n}} \left[ (q_{i+k}-q_{i+1}) \sin \frac{\pi (k+1)}{n}  -  (q_{i+k+1}-q_i) \sin \frac{\pi (k-1)}{n}  \right]&.
\end{split}
\end{equation*}
We want to show that the sequence 
\begin{equation} \label{test}
c_i := (q_{i+k}-q_{i+1}) \sin \frac{\pi (k+1)}{n}  -  (q_{i+k+1}-q_i) \sin \frac{\pi (k-1)}{n}
\end{equation}
must change sign at least four times.

First, observe that $c_i$ is $\ell_2$-orthogonal to the constant sequence $(1,\ldots,1)$, that is, $\sum_{i=1}^n  c_i =0$; hence $c_i$ must have sign changes.

Next, I claim that $c_i$ is $\ell_2$-orthogonal to the  sequence $\sin (2\pi i/n)$.
Indeed 
\begin{equation*}
\begin{split}
\sum_{i=1}^{n} c_i \sin \frac{2\pi i}{n} = \sum_{i=1}^{n} &q_i \sin \frac{\pi (k+1)}{n}
\left( \sin \frac{2\pi (i-k)}{n} - \sin \frac{2\pi (i-1)}{n} \right)\\
- &q_i  \sin \frac{\pi (k-1)}{n} \left( \sin \frac{2\pi (i-k-1)}{n}  + \sin \frac{2\pi i}{n} \right).
\end{split}
\end{equation*}
Hence twice the coefficient of $q_i$ on the right hand side equals
\begin{equation*}
\begin{split}
&\sin \frac{\pi (k+1)}{n} \sin \frac{\pi (1-k)}{n} \cos \frac{\pi (2i-k-1)}{n}\\ 
&+ \sin \frac{\pi (k-1)}{n} \sin \frac{\pi (1+k)}{n} \cos \frac{\pi (2i-k-1)}{n} =0,
\end{split}
\end{equation*}
as needed.

Similarly, $c_i$ is $\ell_2$-orthogonal to the  sequence $\cos (2\pi i/n)$.

Finally, if $c_i$ changes sign only twice, one can find a linear combination 
$$
c+ a \sin \frac{2\pi i}{n} + b \cos \frac{2\pi i}{n},
$$
a discrete first harmonic, that changes sign at the same positions as $c_i$. This ``first harmonic" has no other sign changes, so its signs coincide with those of $c_i$. But it is also orthogonal to $c_i$, a contradiction.
\proofend

\section{Back to four vertices, and another problem}
Perhaps the oldest result in the spirit of  the four vertex-like theorem is the Legendre-Cauchy Lemma (which is about 100 years older than the theorem of Mukhopadhyaya):
{\it if two convex polygons in the plane have equal respective side length, then the cyclic sequence of the differences of their respective angles has at least four sign changes}. 

A  version of this lemma in spherical geometry is the main ingredient of the proof of the Cauchy rigidity theorem ({\it convex polytopes with congruent corresponding faces are congruent to each other}); interestingly,  its original proof contained an error that remained unnoticed for nearly a century, see, e.g., chapters 22 and 26 of \cite{Pak}.

The values of the angles in the formulation of the Legendre-Cauchy Lemma can be replaces by the lengths of the short, skip-a-vertex, diagonals of the respective polygons: with fixed side lengths, the angles depend monotonically on these diagonals.

In particular, one may assume that the polygons are equilateral, e.g., each side has unit length. In this formulation, the Legendre-Cauchy Lemma becomes an analog of the $k=1$ case of  Theorem above, with the determinants $a_i=\det (V_{i-1},V_{i+1})$  replaced by the lengths $|V_{i+1}-V_{i-1}|$. This prompts to ask another question.
\medskip

\noindent {\bf Problem 2.}
{\it 
Given two equilateral convex $n$-gons, for which $k$ must the cyclic sequence $|V_{i+k}-V_{i-1}|$ have at least four sign changes?
}

\bigskip
{\bf Acknowledgements}. It is a pleasure to acknowledge stimulating discussions with S. Morier-Genoud, V. Ovsienko, I. Pak, and R. Schwartz. 
Many thanks to M. Cuntz for providing his (counter)-examples. This work was supported by NSF grant DMS-1510055.

\end{document}